# STABILIZABILITY AND PERCOLATION IN THE INFINITE VOLUME SANDPILE MODEL

By Anne Fey, Ronald Meester and Frank Redig

*VU—University Amsterdam, VU—University Amsterdam and Universiteit Leiden*

We study the sandpile model in infinite volume on $\mathbb{Z}^d$. In particular, we are interested in the question whether or not initial configurations, chosen according to a stationary measure $\mu$, are $\mu$-almost surely stabilizable. We prove that stabilizability does not depend on the particular procedure of stabilization we adopt. In $d = 1$ and $\mu$ a product measure with density $\rho = 1$ (the known critical value for stabilizability in $d = 1$) with a positive density of empty sites, we prove that $\mu$ is not stabilizable.

Furthermore, we study, for values of $\rho$ such that $\mu$ is stabilizable, percolation of toppled sites. We find that for $\rho > 0$ small enough, there is a subcritical regime where the distribution of a cluster of toppled sites has an exponential tail, as is the case in the subcritical regime for ordinary percolation.

**1. Introduction.** The sandpile model was originally introduced as a dynamical model to illustrate the concept of self-organized criticality [1]. The model is defined on a finite subset $\Lambda$ of $\mathbb{Z}^d$, in discrete time. It starts with a stable configuration, that is, every site has a nonnegative height of at most $2d - 1$ sand grains. Every discrete time step, an addition of one sand grain is made to a random site. If this site becomes *unstable*, that is, has at least $2d$ grains, it *topples*, that is, it gives one grain to each neighbor. This may cause other sites to become unstable, and the topplings continue until every site is stable again. The total of all necessary topplings is called an avalanche, so after the avalanche we have reached the new configuration. This is possible in a finite number of topplings because, at the boundary of $\Lambda$, grains are dissipated. This model is abelian: the obtained configuration is independent of the order of topplings.

This sandpile model is said to exhibit self-organized critical behavior, for the following reasons. As the model evolves in time, it reaches a stationary state that is characterized, in the large-volume limit, by long-range height correlations and power law statistics for avalanche sizes, and thus reminds one of critical behavior in statistical mechanical models. However, the sandpile model evolves naturally toward this critical state, without apparent tuning of any parameters.

This seeming contrast has been discussed in [3, 4, 8]; it is argued that the model definition in fact does involve tuning. Namely, the instantaneousness of topplings, and the vanishing of dissipation as $\Lambda \uparrow \mathbb{Z}^d$, can be viewed as a tuning of the addition and dissipation rate to 0 respectively. This tuning then would ensure that the model evolves toward the critical point of a parametrized, nondynamical sandpile model, which can informally be described as follows. We start with an initial height configuration on $\mathbb{Z}^d$ (not necessarily stable) according to a translation invariant probability measure with density $\rho$, which is the expected height, or number of sand grains per site. We keep toppling until there are no more unstable sites. If this is possible with a finite number of topplings per site, then we obtain the final configuration, and the initial configuration is said to be *stabilizable*. This version of the sandpile model was introduced in [3], and mathematically investigated in [4, 8]. Results so far obtained are as follows: for $d = 1$, any translation invariant probability measure with density $\rho < 1$ is stabilizable, any translation invariant probability measure with density $\rho > 1$ is not stabilizable and, for $\rho = 1$, there are cases of stabilizability and nonstabilizability; see [8]. For general $d$, any translation invariant probability measure with density $\rho < d$ is stabilizable, and any translation invariant probability measure with density $\rho > 2d - 1$ is not stabilizable, in between $d$ and $2d - 1$ there are nonstabilizable and stabilizable cases [4].

The present paper continues this investigation and from here on, when we talk about the sandpile model, we mean the version in infinite volume. In Section 2 we introduce notation, introduce general toppling procedures and discuss stabilizability issues. In this section we also prove that if a random initial configuration is stabilizable, then the expected height (density) is conserved by stabilization. In Section 3 we define critical values, and investigate the behavior at the critical point in $d = 1$. We find that configurations chosen according to a nondegenerate product measure are a.s. not stabilizable. In Section 4 we investigate phase transitions for the sandpile model from a new viewpoint: we consider, for stabilizable configurations, percolation of the collection of toppled sites. We look for a critical $\rho$, not necessarily equal to $\rho_c$ mentioned above, such that, for all $\rho$ below this value, there is no infinite cluster of toppled sites.

For a general class of initial distributions and $\rho$ small enough, we find a subcritical regime where not only there is a.s. no infinite cluster of toppled sites, but the distribution of the cluster size has an exponential tail.

This corresponds to the subcritical regime for ordinary percolation, thus strengthening the idea of a critical phase transition.

**2. Toppling procedures and stabilizability.** Denote by $\mathcal{X} = \mathbb{N}^{\mathbb{Z}^d}$ the set of all height configurations and by $\Omega = \{0, 1, \ldots, 2d-1\}^{\mathbb{Z}^d}$ the set of stable height configurations.

A toppling at site $x$ applied to the configuration $\eta \in \mathcal{X}$ is denoted by $\theta_x(\eta)$ and defined via

$$\theta_x(\eta)(y) = \begin{cases} \eta(y) - 2d, & \text{if } y = x, \\ \eta(y) + 1, & \text{if } |y - x| = 1, \\ \eta(y), & \text{otherwise}. \end{cases} \tag{1}$$

A toppling at site $x \in \mathbb{Z}^d$ is called *legal* for the configuration $\eta \in \mathcal{X}$, if it is applied to an unstable site, that is, if $\eta(x) \geq 2d$.

The above definition of a toppling gives rise to the definition of the *toppling matrix* $\Delta$ associated to the sandpile model. This is a matrix indexed by sites $x, y \in \mathbb{Z}^d$, with entries

$$\Delta_{x,y} = 2d \mathbf{1}_{x=y} - \mathbf{1}_{|x-y|=1}.$$

With this definition, and with $\delta_x$ defined to be the vector with entry 1 at $x$ and entry 0 in all other positions, we can write

$$\theta_x(\eta) = \eta - \Delta \delta_x.$$

DEFINITION 2.1. A *toppling procedure* is a measurable map (with respect to the usual Borel sigma-algebra's)

$$T : [0, \infty) \times \mathbb{Z}^d \times \mathcal{X} \to \mathbb{N} \tag{2}$$

such that, for all $\eta \in \mathcal{X}$,

(a) for all $x \in \mathbb{Z}^d$,

$$T(0, x, \eta) = 0.$$

(b) for all $x \in \mathbb{Z}^d$,

$$t \mapsto T(t, x, \eta)$$

is right-continuous and nondecreasing with jumps of size at most one, that is, for all $t > 0$, $x \in \mathbb{Z}^d$, $\eta \in \mathcal{X}$, we have

$$T(t, x, \eta) - T(t-, x, \eta) \leq 1.$$

(c) for all $x \in \mathbb{Z}^d$, in every finite time interval, there are finitely many jumps at $x$.

(d) $T$ does not contain an "infinite backward chain of topplings," that is, there is no infinite chain of topplings at sites $x_i$, $i = 1, 2, \ldots$, occurring at times $t_i > t_{i+1} > \cdots$, where for all $i$, $x_{i+1}$ is a neighbor of $x_i$.

Note that condition (d) is only relevant in continuous time. We interpret $T(t, x, \eta)$ as the number of topplings at site $x$ in the time interval $[0, t]$, when $T$ is applied to the initial configuration $\eta \in \mathcal{X}$. The vector of all such numbers at time $t$ is denoted by $T(t, \cdot, \eta)$. We say that for all $t$ such that $T(t-, x, \eta) < T(t, x, \eta)$, site $x$ topples at time $t$.

If $T$ is a toppling procedure, then for $\eta \in \mathcal{X}$, $t > 0$, we call
$$\Theta_t^\eta(T) = \{x \in \mathbb{Z}^d : T(t, x, \eta) > T(t-, x, \eta)\}$$
the set of sites that topple at time $t > 0$ (for initial configuration $\eta$).

DEFINITION 2.2. Let $T$ be a toppling procedure. The configuration $\eta_t$ at time $t > 0$ associated to $T$ and initial configuration $\eta \in \mathcal{X}$ is defined to be

$$\eta_t = \eta - \Delta T(t, \cdot, \eta). \tag{3}$$

DEFINITION 2.3. A toppling procedure $T$ is called *legal* if for all $\eta \in \mathcal{X}$, for all $t > 0$ and for all $x \in \Theta_t^\eta(T)$, $\eta_{t-}(x) \geq 2d$.

In words, this means that in a legal toppling procedure, only unstable sites are toppled.

DEFINITION 2.4. (a) A toppling procedure $T$ is called *finite* for initial configuration $\eta \in \mathcal{X}$, if for all $x \in \mathbb{Z}^d$,

$$T(\infty, x, \eta) := \lim_{t \to \infty} T(t, x, \eta) = \sup_{t \geq 0} T(t, x, \eta) \tag{4}$$

is finite.

(b) A legal toppling procedure $T$ is called *stabilizing* for initial configuration $\eta \in \mathcal{X}$ if it is finite and if the limit configuration $\eta_\infty$, defined by

$$\eta_\infty = \eta - \Delta T(\infty, \cdot, \eta), \tag{5}$$

is stable.

A *random toppling procedure* is a random variable with values in the set of toppling procedures. This can also be viewed as a measurable map

$$T : [0, \infty) \times \mathbb{Z}^d \times \mathcal{X} \times \hat{\Omega} \to \mathbb{N},$$

where $\hat{\Omega}$ denotes a probability space, and such that for all $\omega \in \hat{\Omega}$ except a set of measure 0, $T(\cdot, \cdot, \cdot, \omega)$ is a toppling procedure.

DEFINITION 2.5. A toppling procedure is called *stationary* if for all $t$, the distribution of $T(t, \cdot, \eta)$ is translation invariant when we choose $\eta$ according to a translation invariant probability measure.

Next we discuss some examples. These examples have in common that for every $t$, if $\eta_t$ contains unstable sites, then each of these sites will topple within finite time almost surely. As a consequence, for every $\eta$, these toppling procedures are either stabilizing or infinite. Moreover, if they are infinite, then $T(\infty, x, \eta) = \infty$ for every $x$. This can be seen as follows: if there is one site $x$ that topples infinitely many times, then the neighbors of $x$ receive infinitely many sand grains. Therefore, these neighbors need to topple infinitely many times, etc.

1. *Markov toppling processes.* These are examples of random stationary toppling procedures and are defined as follows. Each site $x \in \mathbb{Z}^d$ has a Poisson clock (different clocks are independent) with rate one. When the clock at site $x$ rings at time $t$ and in the configuration $\eta_{t-}$, $x$ is unstable, then $x$ is toppled. More formally, the configuration $\eta_t$ of (3) is evolving according to the Markov process with generator, defined on local functions $f : \mathcal{X} \to \mathbb{R}$ via

$$Lf(\eta) = \sum_x \mathbf{1}_{\eta(x) \geq 2d}(f(\theta_x \eta) - f(\eta)).$$

   It is not hard to see that this procedure satisfies all requirements, in particular, (d) of Definition 2.1.

   It is also possible to adapt the rate at which unstable sites are toppled according to their height. In that case the Markov process becomes

$$L_c f(\eta) = \sum_x \mathbf{1}_{\eta(x) \geq 2d} c(\eta(x))(f(\theta_x \eta) - f(\eta)),$$

   where $c : \mathbb{N} \to \mathbb{R}$ has to satisfy certain conditions in order to make the process well-defined.

2. *Toppling in nested volumes.* This is a deterministic, discrete time toppling procedure. Choose a sequence $V_n \subset V_{n+1} \subset \mathbb{Z}^d$ such that $\bigcup_n V_n = \mathbb{Z}^d$, but all $V_n$ contain finitely many sites. We start toppling all the unstable sites in $V_0$ until the configuration in $V_0$ has no unstable sites left, then we do the same with $V_1$, etc. We put this into the framework of Definition 2.1 as follows. At time $t = 1$, we topple all the unstable sites in $V_0$ once, at time $t = 2$, we topple all the unstable sites in $V_0$ if there are still unstable sites left after the topplings at time $t = 1$, etc., until at time $t = t(V_0, \eta)$, no unstable sites are left in $V_0$; we then start toppling at time $t = t(V_0, \eta)$ all unstable sites in $V_1$, etc. Since the volumes $V_n$ are finite, all $t(V_n, \eta)$ are finite.

We will use this procedure several times, but for ease of notation we will reparametrize time such that $V_n$ is stabilized at time $n$ instead of at time $t(V_n, \eta)$.

3. *Topplings in parallel.* Topplings in parallel consists simply in toppling at time $t$ all unstable sites of $\eta_{t-1}$ once. This toppling procedure is discrete time, deterministic and stationary.

4. *Topplings in waves.* This procedure is only used for initial configurations having a single unstable site, say, at $x \in \mathbb{Z}^d$. For the formal definition of the toppling function $T(t, x, \eta)$, we put it equal to zero for configurations $\eta$ that contain more than one unstable site. Toppling in waves is defined for the sandpile model on a finite grid as follows [6]: at $t = 1$, we topple $x$ once and, subsequently, all other sites that become unstable. All these topplings form the first wave. If after these topplings, $x$ is still unstable, then at $t = 2$ we perform the second wave, etc. In each wave, no site topples more than once.

This does not fit into our framework, because in each wave all topplings, except the toppling at $x$, are illegal. Nevertheless, we want each wave to be completed in finite time. Therefore, we define topplings in waves as follows: At $t = 1$, we topple site $x$ once. Then for $i = 2, 3, \ldots$, at times $2 - \frac{1}{i}$, we consecutively topple all sites that are unstable except site $x$. That way the first wave is completed at time $t = 2$. All other waves proceed similarly. Since in each wave no site topples more than once, this procedure is well defined.

DEFINITION 2.6. (a) A configuration $\eta \in \mathcal{X}$ is called *stabilizable* if there exists a stabilizing legal toppling procedure.

(b) A probability measure $\mu$ on $(\mathcal{X}, \mathcal{F})$ is called stabilizable if $\mu$-almost every $\eta$ is stabilizable.

EXAMPLE 2.7. We give an example of a configuration that is not stabilizable: Consider the configuration $\xi$ in $\mathbb{Z}$ where all sites have height 1, except the origin, which has height 2. From trying out by hand, it should become clear that this configuration is not stabilizable. We may choose to topple in waves, since there is only one unstable site. In our case, in each wave every site topples exactly once, so that after each wave we obtain the same configuration; hence, there are infinitely many waves. Alternatively, we may choose to topple in parallel. Then in our case, the height of the origin alternates between 0 and 2, so that the origin topples infinitely many times. From the forthcoming Theorem 2.8, we can use either (3) or (4) to conclude that $\xi$ is not stabilizable.

In [4], Definition 2.4, stabilizability is defined in terms of toppling in nested volumes, and in [4], Lemma 6.12, it is proved that this definition of stabilizability is equivalent for this toppling procedure and the Markov toppling procedure. Here, we extend this result: we prove that if $\eta$ is stabilizable, then irrespective of what legal toppling procedure we choose, it will always be finite, and irrespective of what stabilizing procedure we choose, we always obtain the same stable configuration. On the other hand, if we find one infinite legal toppling procedure for $\eta$, then we know that $\eta$ is not stabilizable. This can also be concluded from the existence of a legal toppling procedure in which every site topples at least once.

Let $T$, $T'$ be two toppling procedures, which are finite for initial configuration $\eta$. Then we write

$$T' \preceq_\eta T,$$

if for all $x \in \mathbb{Z}^d$

$$T'(\infty, x, \eta) \leq T(\infty, x, \eta).$$

THEOREM 2.8. *Let $T$, $T'$ be two legal toppling procedures, which are both finite for initial configuration $\eta$:*

1. *If $T$ is stabilizing for $\eta$, then*

(6) $$T' \preceq_\eta T.$$

2. *If $T$ and $T'$ are two stabilizing toppling procedures for $\eta$, then for all $x \in \mathbb{Z}^d$,*

$$T'(\infty, x, \eta) = T(\infty, x, \eta).$$

   *In particular, this means that for stabilizable $\eta$, the limit configuration $\eta_\infty$ is well defined.*
3. *For stabilizable $\eta \in \mathcal{X}$, there does not exist a nonfinite legal toppling procedure.*
4. *If $T$ is stabilizing for $\eta$, then there is at least one site $x$ that does not topple, that is, there is at least one site $x$ for which $T(\infty, x, \eta) = 0$.*

PROOF. The proof of Statement 1 is inspired by an argument that appears in [2] and [9] in the context of finite grids or discrete time toppling procedures.

For every $x$, we define a time $\tau_x := \sup\{t : T'(t, x, \eta) \leq T(\infty, x, \eta)\}$, and we call all topplings in $T'$ that occur at times strictly larger than $\tau_x$ "extra" topplings. We suppose the converse of Statement 1, that is, we suppose that there is at least one extra toppling.

Suppose an extra toppling occurs at site $y$, at time $t_y < \infty$. Then just before time $t_y$, the number of topplings at site $y$ is at least $T(\infty, y, \eta)$.

Moreover, in order for this extra toppling to be legal, site $y$ must be unstable just before time $t_y$. Thus, we find, following $T'$, that

$$2d \leq \eta_{t_y-}(y) = \eta(y) - (\Delta T'(t_y-,\cdot,\eta))(y)$$
$$= \eta(y) - 2dT'(t_y-,y,\eta) + \sum_{x \sim y} T'(t_y-,x,\eta)$$
$$\leq \eta(y) - 2dT(\infty,y,\eta) + \sum_{x \sim y} T'(t_y-,x,\eta),$$

where the sum $\sum_{x \sim y}$ runs over all neighbors of $y$. Since $T$ is stabilizing, we have

$$2d > \eta(y) - 2dT(\infty,y,\eta) + \sum_{x \sim y} T(\infty,x,\eta),$$

so for at least one $x \sim y$, $T'(t_y-,x,\eta) > T(\infty,x,\eta)$. In other words, for an extra toppling at site $s$ to be legal, it is necessary that it is preceded by at least one extra toppling at one of its neighbors. Then for this extra toppling, we can make the same observation. Continuing this reasoning, we find that in order for the extra toppling at $s$ to be legal, we need an infinite backward chain of extra topplings, occurring in finite time. But then $T'$ does not satisfy item (d) of Definition 2.1. This proves Statement 1.

To prove Statement 2, we simply observe that if $T$ and $T'$ are both stabilizing, then according to the above, $T' \preceq_\eta T$ and $T \preceq_\eta T'$, so that they must be equal.

To prove Statement 3, let $T$ be a stabilizing toppling procedure, and $T''$ a nonfinite legal toppling procedure. Since $T''$ is nonfinite, there exists $x \in \mathbb{Z}^d$ such that $T''(t,x,\eta) \uparrow \infty$ as $t \uparrow \infty$. For some $w < \infty$, we define $T''_w$ as follows: for all $t \leq w$, $T''_w(t,\cdot,\eta) = T''(t,\cdot,\eta)$, but for all $t > w$, $T''_w(t,\cdot,\eta) = T''(w,\cdot,\eta)$. In words, $T''_w$ performs all topplings according to $T''$ up to time $w$, but then stops toppling. $T''_w$ is a finite legal toppling procedure by item (c) of Definition 2.1 and, hence, by Statement 1 of this theorem, $T''_w(\infty,x,\eta) \leq T(\infty,x,\eta)$. By letting $w \to \infty$, we obtain a contradiction.

To prove Statement 4, suppose that there is a stabilizing toppling procedure $T$ such that $T(\infty,x,\eta) > 0$ for all $x$. For every $x$, we call the toppling that occurs according to $T$ at time $t_x := \min\{t : T(t,x,\eta) = T(\infty,x,\eta)\}$ the "last" toppling. Since $T$ is stabilizing, $t_x$ is finite for all $x$.

We define $\bar{T}$ as

$$\bar{T}(t,x,\eta) := \min\{T(t,x,\eta), T(\infty,x,\eta) - 1\},$$

so that for all $x$, $\bar{T}(\infty,x,\eta) = T(\infty,x,\eta) - 1$. In words, $\bar{T}$ contains all topplings according to $T$ except the last one at each site. Note that $\bar{T}$ is a finite, but not a priori legal toppling procedure. However, we have

$$\eta - \Delta\bar{T}(\infty,\cdot,\eta) = \eta - \Delta T(\infty,\cdot,\eta) = \eta_\infty,$$

so that after all topplings according to $\bar{T}$, we have a stable configuration. Now the argument proceeds as in the proof of Statement 1: we have, for some site $v$,

$$2d \leq \eta_{t_v-}(v) = \eta(v) - 2dT(t_v-, v, \eta) + \sum_{x \sim v} T(t_v-, x, \eta)$$
$$= \eta(v) - 2d\bar{T}(\infty, v, \eta) + \sum_{x \sim v} T(t_v-, x, \eta),$$

whereas

$$2d > \eta(v) - 2d\bar{T}(\infty, v, \eta) + \sum_{x \sim v} \bar{T}(\infty, x, \eta).$$

Similarly, as in the proof of Statement 1, we conclude that for the last toppling at $v$ to occur legally, it must have been preceded by an infinite backward chain of last topplings, so that $T$ cannot satisfy item (d) of Definition 2.1. □

REMARK 2.9. Note that if $\mu$ is stabilizable and ergodic, then the induced measure on limit configurations is also ergodic since it is a factor of $\mu$.

We now prove that a finite legal toppling procedure conserves the density. From here on, we will denote by $\mathbb{E}_\mu$, $\mathbb{P}_\mu$ expectation resp. probability with respect to $\mu$.

LEMMA 2.10. *Let $\mu$ be a translation invariant probability measure on $\mathcal{X}$ such that $\mathbb{E}_\mu(\eta(0)) = \rho < \infty$. Suppose, furthermore, that $\mu$ is stabilizable. Then the expected height is conserved by stabilization, that is,*

$$\mathbb{E}_\mu(\eta_\infty(0)) = \rho.$$

*Moreover, if $\mu$ is a translation invariant probability measure on $\mathcal{X}$ such that $\mathbb{E}_\mu(\eta(0)) = \infty$, then $\mu$ is not stabilizable.*

PROOF. Using ergodic decomposition, we can assume that $\mu$ is ergodic. We start with the case $\mathbb{E}_\mu(\eta(0)) = \rho < \infty$. Without loss of generality, we assume that the toppling procedure that stabilizes $\mu$ is stationary, and, is moreover, such that for all $t$, $\mathbb{E}_\mu(T(t, x, \eta)) < \infty$ [we can, e.g., choose the Markov toppling procedure, where $T(t, x, \eta)$ is dominated by a Poisson process].

At time $t$, we then have

$$\eta_t(x) = \eta(x) - \sum_y \Delta_{x,y} T(t, y, \eta),$$

which upon integrating over the distribution of $\eta$ gives

$$\mathbb{E}_\mu(\eta_t(x)) = \mathbb{E}_\mu(\eta(x)) = \rho.$$

Therefore, using Fatou's lemma,

$$\rho_\infty := \mathbb{E}_\mu(\eta_\infty(0)) = \mathbb{E}_\mu\left(\lim_{t\to\infty}\eta_t(0)\right) \leq \liminf_{t\to\infty}\mathbb{E}_\mu(\eta_t(0)) = \rho.$$

The inequality $\mathbb{E}_\mu(\eta_\infty(0)) \geq \rho$ is proved in [4]; we give a somewhat different argument here. Let $X_n$ denote the position of a simple random walk starting at the origin (independent of $\eta$), and denote $\mathbb{E}_{rw}, \mathbb{P}_{rw}$ expectation and probability with respect to this random walk. We start by choosing a stabilizable $\eta$ with limit $\eta_\infty$, and for a moment we consider this $\eta$ and $\eta_\infty$ fixed.

From the relation

$$\eta_\infty(x) = \eta(x) - \sum_y \Delta_{x,y} T(\infty, y, \eta),$$

we obtain

$$(7) \quad \frac{1}{n}\mathbb{E}_{rw}\left(\sum_{k=0}^{n-1}(\eta_\infty(X_k) - \eta(X_k))\right) = \frac{2d}{n}\mathbb{E}_{rw}(T(\infty, X_n, \eta)) - \frac{2d}{n}T(\infty, 0, \eta).$$

By letting $n \to \infty$, using that $T(\infty, 0, \eta) < \infty$ by assumption, and that $\mathbb{E}_{rw}(T(\infty, X_n, \eta)) \geq 0$, this leads to

$$(8) \quad \liminf_{n\to\infty}\frac{1}{n}\mathbb{E}_{rw}\left(\sum_{k=0}^{n-1}\eta_\infty(X_k)\right) \geq \liminf_{n\to\infty}\frac{1}{n}\mathbb{E}_{rw}\left(\sum_{k=0}^{n-1}\eta(X_k)\right).$$

If we now finally choose $\eta$ according to $\mu$, which is ergodic, then the limiting measure is also ergodic according to Remark 2.9. By ergodicity of the scenery process $\{\eta(X_n): n \in \mathbb{N}\}$ (see, e.g., [4], Proposition 8.1), it follows that for $\mu$-a.e. $\eta$, the right-hand side is equal to $\rho$, and the left-hand side is equal to $\rho_\infty$. This proves that $\rho_\infty \geq \rho$.

Finally, if $\mathbb{E}_\mu(\eta(0)) = \infty$, then the right-hand side of (8) goes to $+\infty$ as $n \to \infty$. Therefore, if $\mu$ would be stabilizable, then $\eta_\infty$ would have infinite density, a contradiction. □

**3. Criticality and critical behavior.** Let $\mathcal{P}(\mathcal{X})$ denote the set of all translation invariant probability measures on $(\mathcal{X}, \mathcal{F})$. We say that a subset $\mathcal{M}$ of $\mathcal{P}(\mathcal{X})$ is *density complete* if for all $\rho \in [0, \infty)$ there exists $\mu \in \mathcal{M}$ such that $\mu(\eta(0)) = \rho$.

Let $\mathcal{M} \subseteq \mathcal{P}(\mathcal{X})$ be density complete. We define the $\mathcal{M}$-*critical density* for stabilizability to be

$$(9) \quad \rho_c(\mathcal{M}) = \sup\{\rho > 0 : \forall \mu \in \mathcal{M} \text{ with } \mu(\eta(0)) = \rho, \mu \text{ is stabilizable}\}.$$

Of course, it can be questioned whether the density is the only relevant parameter distinguishing between stabilizability and nonstabilizability. It is certainly the most natural parameter, and is considered in the numerical experiments of [3]. In [4] a related notion of *maximal stabilizability* is introduced.

It is clear that $\rho_c(\mathcal{M}) \leq \rho_c(\mathcal{M}')$ for $\mathcal{M} \supseteq \mathcal{M}'$. Natural choices for $\mathcal{M}$ are a one-parameter family of product measures such as the set of Poisson product measures with parameter $\rho$, the set of all product measures or simply $\mathcal{M} = \mathcal{P}(\mathcal{X})$.

The following results are reformulations of results in [4] and [8].

THEOREM 3.1.  (a) *For* $\mathcal{M} = \mathcal{P}(\mathcal{X})$, *and for all* $d$,
$$\rho_c(\mathcal{M}) = d.$$
(b) *For all* $\mathcal{M}$ *density complete, we have*
$$d \leq \rho_c(\mathcal{M}) \leq 2d - 1.$$
*In particular, when* $d = 1$ *and for all* $\mathcal{M}$ *density complete, we have*
$$\rho_c(\mathcal{M}) = 1.$$

We now specialize to the case $d = 1$. Accordingly, let $\mu$ be a one-dimensional translation invariant product measure with density $\rho = 1$. From Theorem 3.1, we know that for $\rho < 1$, $\mu$ is stabilizable, and that for $\rho > 1$, it is not. The next result deals with the critical case $\rho = 1$.

THEOREM 3.2. *Let* $\mu$ *be a one-dimensional product measure with* $\rho = 1$ *such that* $\mu(\eta(0) = 0) > 0$. *Then* $\mu$ *is not stabilizable.*

Our strategy will be to show that there a.s. exists a nonfinite legal toppling procedure. This implies, using Theorem 2.8, that $\mu$ is not stabilizable. In order to do so, we will use topplings in nested volumes, but for the proof it will be important to define an intermediate toppling procedure, during which we only stabilize in volumes of the form $[0, n]$, that is, we increase the stabilized volume only to one side. After stabilization of the interval $[0, n]$, the outer boundary sites $-1$ and $n + 1$ possibly contain, on top of their original height, extra grains that were removed from the interval $[0, n]$ during stabilization; all other sites outside $[-1, n+1]$ still have their original height.

For a while we concentrate on this one-sided procedure. In this section we denote by $\eta_n$ the configuration that results from stabilizing in the interval $[0, n]$. As in the nested volumes toppling procedure, we re-define *time* as to match this notation: stabilization of $[0, n]$ takes place at time $n$ so that $\eta_n$ is the configuration reached at time $n$.

We will work with the number and positions of empty sites of $\eta_n$ in $[0, n]$, and we will call such an empty site "a 0" of $\eta_n$. In Figure 1 we illustrate the dynamics of the 0's in this procedure. Time (in the new sense) is plotted vertically going downward; space is plotted horizontally. At every time, when you look horizontally, the black dots you see represent the positions of the 0's at that time in the interval $[0, n]$. In addition, the outer boundary sites of $[0, n]$ are also colored black. The configuration outside the stabilized interval is not shown. Thus, the picture does not give complete information about the configuration $\eta_n$, it only shows the positions of the 0's and the width of the stabilized interval. One can clearly follow the trajectory of the leftmost 0. For instance, just before time 4000, the leftmost 0 starts to move, and reaches the outer boundary somewhere between times 6000 and 7000. Shortly after that, a new leftmost 0 starts to move to the right, etc.

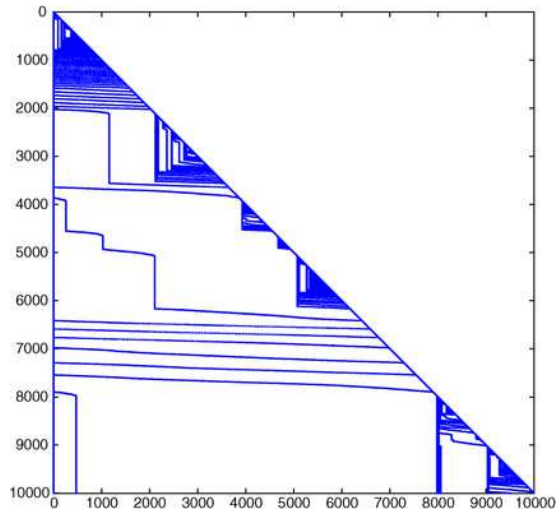

Fig. 1. *The first 10,000 time steps for stabilizing $\eta$ according to* Poisson(1) *product measure in nested volumes $[0, n]$. The y-axis represents time. A black dot indicates an empty site (a "0"). In addition, the outer boundary sites of $[0, n]$ are also colored black. See the text for further explanation.*

Our strategy is to show that during the one-sided procedure, despite the fact that infinitely often new 0's are created, we infinitely often encounter a configuration that does not contain a 0. Every time this occurs, there is a fixed positive lower bound for the probability that the origin topples. This will then imply that the origin topples infinitely many times a.s.

In order to show that infinitely often there are no 0's, we need to analyze the dynamics of the 0's in great detail. We are going to view the 0's as objects that can move, disappear or be created. In order to precisely define these terms, we organize the topplings that occur in time step $n$ into waves. If at time $n - 1$ site $n$ is unstable (at time $n - 1$, all sites $0, \ldots, n - 1$ are stable), then in the first wave at time $n$, we topple site $n$ once and then all other sites in $[0, n]$ that become unstable, except site $n$ again. If after this wave site $n$ is still unstable, then the second wave starts, etc. We will number the waves $k = 0, \ldots, K$, and call $\tilde{\eta}_{n-1,k}$ the configuration after the $k$th wave, so that $\tilde{\eta}_{n-1,0} = \eta_{n-1}$ and $\tilde{\eta}_{n-1,K} = \eta_n$. Depending on the position of the rightmost 0 after wave $k - 1$, wave $k$ has the following effect:

1. If the rightmost 0 of $\tilde{\eta}_{n-1,k-1}$ is at site $n - 1$, and after wave $k$ site $n$ is not empty, then the number of 0's has decreased by 1;
2. If there are no 0's in $\tilde{\eta}_{n-1,k-1}$, then all sites in $[0, n]$ topple, after which there is a 0 at the origin, and site $-1$ has gained a grain.
3. In all other cases, if the rightmost 0 is at position $x - 1$ (so that $x$ is the leftmost site that topples during wave $k$), then site $x - 1$ gains a grain and site $x$ loses one. In addition, site $n$ loses one grain and site $n + 1$ gains one.

These observations inspire the following definition.

DEFINITION 3.3. Let, at time $n$, $K$ be the number of waves. If $K > 0$, then let, in wave $k$, $x$ be the leftmost site that topples:

- If $x = n$, and after wave $k$ site $n$ is not empty, we say that the 0 at site $n - 1$ *disappears*.
- If $x = 0$, we say that a new 0 is *created at the origin*.
- If $x > 0$, and no 0 disappears, we say that the 0 at site $x - 1$ *moves* to site $x$.

If $\eta_{n-1}(n) = 0$ (this implies $K = 0$), we say that a new 0 is *created at the right boundary*.

Since there may be multiple waves in one time step, multiple things can happen to the 0's. However, note that we have the following restrictions: in each wave, only the rightmost 0 can move. For instance, in the example in Figure 1, the 0 that is present at position 472, at time 10,000, has been in that position for almost 2000 time steps, and we cannot be sure if it will ever move again some future time (actually, as our proof will show, it will a.s.). Furthermore, only when, after a previous wave, there are no 0's left can a new 0 be created at the origin. For instance, in the realization in Figure 1, this occurs seven times between $n = 6000$ and $n = 8000$.

We stress that, according to the above definition, we actually identify certain 0's in different time steps. A look at Figure 1 should convince the reader that this is a natural way to view the 0's, even though in order to do so, it is necessary to break the topplings in each time step up into waves to ensure a correct identification. Once a 0 has been created, it exists until it disappears at the right boundary. This may be in the same time step, but it could also require many time steps. During this time, it may move to the right or remain for some periods of time in the same position.

The time intervals between successive instances where the number of 0's is equal to some given number $z$ are not i.i.d. time intervals. However, we will show in the following lemma that, for all $z > 0$, the time intervals, from the moment that, the number of 0's becomes $z + 1$ until the first return to a value that is at most $z$, are i.i.d. time intervals, whose distribution does not depend on $z$. In the proof we use that, for $z > 0$, the number of 0's can only increase from $z$ to $z + 1$ when a new 0 is created at the right boundary. When $z = 0$ we can have that the number of 0's increases because a new 0 is created at the origin, in which case the proof does not apply.

LEMMA 3.4. *Let $z > 0$. Let $Z(n)$ be the number of 0's after time $n$. For $i = 0, 1, \ldots$, let $N_0(z) = 0$, $M_i(z) = \min\{n > N_i(z) : Z(n) \leq z\}$, and $N_i(z) = \min\{n > M_{i-1}(z) : Z(n) = z + 1\}$. Then:*

1. *The random variables $\Delta_i(z) = M_i(z) - N_i(z)$ are i.i.d., for all $i > 0$.*
2. *The distribution of $\Delta_i(z)$ does not depend on $z$, and we denote by $\Delta$ a random variable with this distribution.*
3. *If $\liminf_{n \to \infty} Z(n) < \infty$ a.s., then $\liminf_{n \to \infty} Z(n) \leq 1$ a.s. and $\mathbb{P}(\Delta = \infty) = 0$.*

PROOF. Since at all times $N_i(z)$ a new 0 is created at the right boundary, it must be the case that

$$\eta_{N_i(z)-1}(N_i(z)) = 0.$$

This new 0 will be present until time $M_i(z)$ and, during this time, it cannot move to the left. The key observation is that the dynamics of this new 0 depend only on $\eta(j)$, $j \geq N_i(z)$. In particular, conditioned on the creation of the new 0 at time $N_i(z)$, the value of $M_i(z)$ only depends on these random variables. Since

1. the $M_i(z)$ and $N_i(z)$ are all stopping times, and
2. $\mu$ is a product measure,

it follows that the $\Delta_i(z)$ are i.i.d. random variables for every fixed $z$. Furthermore, it follows immediately that the distribution of $\Delta_i(z)$ is also independent of $z$. This proves the first two parts of the lemma.

We now proceed with Part 3. Suppose that $\liminf_{n \to \infty} Z(n) > 1$ with positive probability, that is, there is a random variable $N$, finite with positive probability, such that $Z(n) > 1$ for all $n > N$. We denote by $N_e(1)$ the total number of time intervals $\Delta_i(1)$ [during which $Z(n) > 1$]. If $N_e(1) < \infty$, then the last interval has infinite length. However, since $\mathbb{P}(\liminf_{n \to \infty} Z(n) > 1) > 0$, it is the case that $\mathbb{P}(N_e(1) < \infty) > 0$. We calculate, using that the $\Delta_i(1)$'s are independent, $\mathbb{P}(N_e(1) < \infty) = \sum_{k=1}^{\infty} \mathbb{P}(N_e(1) = k) = \sum_{k=1}^{\infty} \prod_{i=1}^{k-1} \mathbb{P}(\Delta_i(1) < \infty)\mathbb{P}(\Delta_k(1) = \infty)$. This equals $\sum_{k=1}^{\infty} \mathbb{P}(\Delta < \infty)^{k-1}\mathbb{P}(\Delta = \infty)$, so that we obtain $\mathbb{P}(\Delta = \infty) > 0$.

So far, we showed that $\liminf_{n \to \infty} Z(n) > 1$, with positive probability, implies that $\mathbb{P}(\Delta = \infty) > 0$. Now we show that $\mathbb{P}(\Delta = \infty) > 0$ implies $\liminf_{n \to \infty} Z(n) = \infty$ a.s. Denote by $N_e(z)$ the total number of time intervals $\Delta_i(z)$ [during which $Z(n) > z$], and call $\mathbb{P}(\Delta = \infty) = p$. Similarly to the above computation, we calculate $\mathbb{P}(N_e(z) < \infty) = \sum_{k=1}^{\infty} \mathbb{P}(\Delta < \infty)^{k-1}\mathbb{P}(\Delta = \infty) = \sum_{k=0}^{\infty}(1-p)^k p = 1$. Therefore, $\mathbb{P}(N_e(z) < \infty) = 1$ for all $z > 0$, so that $\liminf_{n \to \infty} Z(n) > z$ a.s. for all $z > 0$. It follows that $\liminf_{n \to \infty} Z(n) = \infty$ a.s. □

PROOF OF THEOREM 3.2. We choose $\eta$ according to $\mu$. As mentioned before, we will show that there a.s. exists a nonfinite legal toppling procedure. We will use toppling in nested volumes $[-m, m]$. However, in order to compare this procedure with the one-sided procedure introduced above, we will reach $\eta_m$ from $\eta$ in the following way: First, we stabilize the interval $[0, m]$. After this step, site $-1$ received a number $A^+(m)$ of grains, and $[0, m]$ contains a number $Z^+(m)$ of 0's. Then, we stabilize the interval $[-m, -2]$, in the same way as we stabilized in $[0, m]$. After this step, site $-1$ received another number $A^-(m)$ of grains, and $[-m, -2]$ contains a number $Z^-(m)$ of 0's. Site $-1$ is now the only possibly unstable site in $[-m, m]$. Finally, we stabilize all of $[-m, m]$. Note that in this description, for every $m$, we obtain $\eta_m$ starting from $\eta$, whereas in the above presented one-sided procedure we obtained $\eta_n$ from $\eta_{n-1}$. The numbers $A^+(m)$ and $A^-(m)$ are nondecreasing in $m$. The sequences $(Z^+(m))$ and $(Z^-(m))$ are independent of each other, and also have the same distribution.

The following discussion will repeatedly involve both the one-sided and the nested volume toppling procedure. To make the distinction clear, we will use indices $n$ or $N$ to refer to time steps for the one-sided procedure, and indices $m$ or $M$ to refer to time steps of the nested volume procedure.

*First case.* We assume that, with positive probability, $\liminf_{n \to \infty} Z^\pm(m) = \infty$. If both liminfs are actually infinite and we apply the right one-sided procedure to $\eta$, then for every $z > 0$ there is a time $N(z)$ such that, for all $n > N(z)$, $[0, n]$ contains at least $z$ 0's. This, however, implies that the leftmost $z - 1$ 0's never move again, which in turn implies that from some $n$ on, grains can never reach site $-1$ again; a similar argument is valid for the left one-sided procedure on the interval $[-n, -2]$. Hence, there is positive probability that both $A^+(m)$ and $A^-(m)$ do not increase anymore eventually, and therefore remain bounded. However, since both $Z^+(m)$ and $Z^-(m)$ tend to infinity, the number of 0's with fixed positions in both the left and right one-sided procedure tends to infinity. This now is incompatible with stabilization, since after toppling site $-1$ in the end, we should (if stabilization occurs) obtain a stable configuration $\eta_\infty$ which should be equal to $\bar{1}$, by Lemma 2.10. However, there are simply not enough grains at $-1$ to fill all the 0's that were created by the one-sided procedures.

*Second case.* We now know that $\liminf_{m \to \infty} Z^\pm(m) < \infty$ a.s. By Lemma 3.4, part 3, we conclude that a.s. $\liminf_{n \to \infty} Z(n) \leq 1$, and $\mathbb{P}(\Delta = \infty) = 0$. This implies that all 0's, possibly except the leftmost one, will eventually disappear. Although the proof of Lemma 3.4, Part 3, does not work when $z = 0$, we can a fortiori conclude that also the leftmost 0 must eventually disappear. Indeed, since all other 0's eventually disappear, and since the occurrence of this event only depends on the configurations to the right of such a 0, it follows that no matter what the configuration to the right of a certain 0 is, it will always disappear eventually. Clearly, this is then also true for the leftmost 0. [Note though that the leftmost 0 may disappear without $Z(n)$ decreasing; if in a time step where the leftmost 0 disappears also the origin topples, then a new 0 is created at the origin.]

Finally, since clearly infinitely many 0's are created at the right boundary, we conclude that infinitely often the leftmost 0 disappears. Now consider one time instant $N'$ such that the leftmost 0 disappears at time $N'$. Whether a new 0 is created at this time depends on the precise value of $\eta_{N'-1}(N')$. Given that this amount is large enough to make the leftmost 0 disappear, we can either have that the origin topples as well, or we can have that the origin does not topple, so that at time $N'$ there are no 0's. In the last case, if $\eta_{N'}(N'+1) \geq 2$, then the origin topples at time $N'+1$. The probability that $\eta_{N'}(N'+1) \geq 2$ is bounded from below by $\mathbb{P}(\eta(N'+1) \geq 2)$. Thus, we have that at every time instant where the leftmost 0 disappears, either the origin topples, or it topples with at least a fixed positive probability one time step later. We conclude that during the one-sided procedure the origin topples infinitely often, so that the procedure is nonfinite. $\square$

If the initial measure $\mu$ satisfies a central limit theorem, then a shorter proof is possible:

THEOREM 3.5. *Let $\mu$ be a translation invariant probability measure on $\mathcal{X}$ such that $\mathbb{E}_\mu(\eta(0)) = 1$ and such that $\frac{1}{\sqrt{n}} \sum_{x=-n}^{n}(\eta(x) - 1)$ converges in distribution, as $n \to \infty$ to a nondegenerate normal random variable. Then $\mu$ is not stabilizable.*

PROOF. Suppose that $\mu$ is stabilizable. Then there exist random variables $N(x)$ such that, for all $x \in \mathbb{Z}$,

$$\eta - \Delta N(x) = \bar{1}, \tag{10}$$

where $\Delta N(x) := \sum_y \Delta_{xy} N(y)$, and where $\overline{1}$ denotes the configuration with height 1 everywhere. Indeed, if $\mu$ is stabilizable, then the only final stable configuration can be the configuration which is constant and equal to 1 [by conservation of density (Lemma 2.10) and stability]. The variables $N(x)$ are then the number of topplings needed at $x$ to stabilize $\eta$. By stationarity of the toppling mechanism, the joint distribution of $N(x)$ is stationary under translations. From (10), we obtain

$$\frac{1}{\sqrt{n}} \sum_{x=-n}^{n} (\eta(x) - 1) = \frac{1}{\sqrt{n}} \left( \sum_{x=-n}^{n} \Delta N(x) \right)$$
(11)
$$= \frac{1}{\sqrt{n}} (N(-n-1) - N(n) + N(n+1) - N(-n)).$$

We claim that the right-hand side of (11) converges to 0 in probability as $n \to \infty$. Indeed, for $\varepsilon > 0$, by stationarity, we have

$$\mathbb{P}_\mu \left( \frac{1}{\sqrt{n}} (N(-n-1) - N(n) + N(n+1) - N(-n)) \geq \varepsilon \right)$$
$$\leq 4 \mathbb{P}_\mu \left( N(0) \geq \frac{\varepsilon \sqrt{n}}{4} \right),$$

which tends to zero by the assumption that $N(0) < \infty$ a.s. This leads to a contradiction since, by assumption, the left-hand side of (11) converges to a nondegenerate normal random variable. $\square$

**4. Sandpile percolation.** We call $\mathcal{T}_t$ the set of all sites that have toppled at least once up to (and including) time $t$, that is, $\mathcal{T}_t = \{x : T_t(x) > 0\}$. Likewise, we introduce the set of nonempty sites at time $t$, $\mathcal{V}_t = \{x : \eta_t(x) > 0\}$, and finally $\mathcal{W}_t = \mathcal{T}_t \cup \mathcal{V}_t$, the set of sites that have toppled or are nonempty at time $t$.

For $\eta$ stabilizable, these sets have a limit, for example, $\mathcal{T}_\infty = \lim_{t \to \infty} \mathcal{T}_t$. We decompose the set $\mathcal{T}_\infty$ in clusters $\mathcal{T}_\infty(x)$, where $\mathcal{T}_\infty(x)$ is the largest connected component of $\mathcal{T}_\infty$ containing $x$. Likewise, we decompose $\mathcal{W}_\infty$ into clusters $\mathcal{W}_\infty(x)$, where $\mathcal{W}_\infty(x)$ is the largest connected component of $\mathcal{W}_\infty$ containing $x$. Sandpile percolation is the study of these clusters.

As in classical percolation, one can define critical densities for the existence or absence of infinite clusters and distinguish between a sub- and supercritical regime. In this section we are interested in the tail of the cluster size distribution

$$\mathbb{P}_\mu(|\mathcal{T}_\infty(0)| \geq n)$$

and in the percolation probability

$$\mathbb{P}_\mu(|\mathcal{T}_\infty(0)| = \infty).$$

For the other sets, definitions and notation are similar. In this section we will need the following large deviation result.

LEMMA 4.1. *Let $\mu_\rho$, $0 < \rho < 1$ be a one-parameter family of translation invariant product measures satisfying the following:*

(i) *For all $t \in \mathbb{R}$, the moment generating function*

(12) $$G_\rho(t) = \mathbb{E}_{\mu_\rho}(e^{t\eta(0)}) < \infty$$

*exists;*

(ii) *The first moment satisfies*

$$\mathbb{E}_{\mu_\rho}(\eta(0)) = \rho;$$

(iii) *If $a_\rho \geq 0$ are nonnegative such that $a_\rho \to \infty$ for $\rho \to 0$, then*

$$G_\rho^{-1}(a_\rho) := \sup\{x : G_\rho(x) \leq a_\rho\}$$

*also tends to $\infty$ as $\rho \to 0$.*

*Then we have, for any sequence $x_1, x_2, \ldots$ of lattice points and any $\varepsilon > 0$,*

(13) $$\limsup_{\rho \to 0} \limsup_{n \to \infty} \frac{1}{n} \log \mathbb{P}_{\mu_\rho}\left(\sum_{i=1}^n \eta(x_i) \geq \varepsilon n\right) = -\infty.$$

PROOF. By the Markov inequality, for any $t \geq 0$, using that under the $\eta(x_i)$ are i.i.d., we obtain

$$\mathbb{P}_{\mu_\rho}\left(\sum_{i=1}^n \eta(x_i) \geq \varepsilon n\right) \leq e^{-\varepsilon n t} G_\rho(t)^n,$$

which gives

(14) $$\frac{1}{n} \log \mathbb{P}_{\mu_\rho}\left(\sum_{i=1}^n \eta(x_i) \geq \varepsilon n\right) \leq -\varepsilon t + \log G_\rho(t).$$

We now show that there exists $t_\rho > 0$ such that

$$\log G_\rho(t_\rho) \leq 1$$

and such that $t_\rho \to \infty$ as $\rho \to 0$. Using that $\eta(0)$ takes only integer values, the elementary inequality $\log(1 + x) \leq x$, the Cauchy–Schwarz inequality and, finally, the Markov inequality, we obtain

(15) $$\begin{aligned}\log G_\rho(t) &\leq \log(1 + \mathbb{E}_{\mu_\rho}(e^{t\eta(0)} \mathbf{1}_{\eta_0 \geq 1})) \\ &\leq \mathbb{E}_{\mu_\rho}(e^{t\eta(0)} \mathbf{1}_{\eta(0) \geq 1}) \\ &\leq (\mathbb{E}_{\mu_\rho}(e^{2t\eta(0)}))^{1/2}(\mathbb{P}_{\mu_\rho}(\eta(0) \geq 1))^{1/2} \\ &\leq G_\rho(2t)^{1/2} \rho^{1/2},\end{aligned}$$

STABILIZABILITY IN SANDPILES 19so we can choose $t_\rho = \frac{1}{2}G_\rho^{-1}(1/\rho^{1/2})$ and, by condition (iii), we have that $t_\rho \to \infty$ as $\rho \to 0$. We can now finish the proof by returning to (14), choosing $t = t_\rho$ in the right-hand side of the inequality, and letting $\rho \to 0$. □

We remark that an elementary computation shows that conditions (i)–(iii) are satisfied when $\mu_\rho$ has Poisson-$\rho$ one-dimensional marginals.

We now first deal with the tail of $|\mathcal{T}_\infty(0)|$.

THEOREM 4.2. (a) *Let $d = 1$ and let $\mu$ be a translation invariant product measure satisfying $\mathbb{E}_\mu(e^{t\eta(0)}) < \infty$ and with density $\rho < 1$. Then there exists a constant $c_1 > 0$ such that*

$$\mathbb{P}_\mu(|\mathcal{T}_\infty(0)| \geq n) \leq e^{-c_1 n}.$$

(b) *Let $d > 1$ and let $\mu_\rho$, $0 < \rho < 1$ be a collection of product measures satisfying conditions* (i)–(iii) *of Lemma 4.1. Then for all $\rho$ sufficiently small, there exists a constant $c_d = c_d(d, \rho) > 0$ such that*

$$\mathbb{P}_{\mu_\rho}(|\mathcal{T}_\infty(0)| \geq n) \leq e^{-c_d n}.$$

For the proof, we need the following result which goes back to at least [7]. For completeness, we give the short proof.

LEMMA 4.3. *Let $\Lambda$ be a subset of $\mathcal{T}_t$, for some toppling procedure. Let $\beta_\Lambda$ be the number of internal bonds in $\Lambda$, that is, bonds with both endpoints in $\Lambda$. Then*

$$\sum_{x \in \Lambda} \eta_t(x) \geq \beta_\Lambda.$$

PROOF. For each internal bond of $\Lambda$, consider the last particle that traversed this bond via a toppling. This particle remains in $\Lambda$ up to time $t$. The result now follows. □

PROOF OF THEOREM 4.2. Let $\mu$ be a product measure in dimension $d \geq 1$. We then have

$$\mathbb{P}_\mu(|\mathcal{T}_\infty(0)| \geq n) = \sum_{m=n}^{\infty} \mathbb{P}_\mu(|\mathcal{T}_\infty(0)| = m) + \mathbb{P}_\mu(|\mathcal{T}_\infty(0)| = \infty).$$

We choose to stabilize in nested boxes $B_k$ of radius $k$. Recall that we reparametrize time so that, at time $k$, the whole box $B_k$ has been stabilized. Then for every $k$, the maximum size of $\mathcal{T}_k(0)$ is $(2k+1)^d$, so that we can rewrite

$$(16) \quad \mathbb{P}_\mu(|\mathcal{T}_\infty(0)| \geq n) = \lim_{k \to \infty} \mathbb{P}_\mu(|\mathcal{T}_k(0)| \geq n) = \lim_{k \to \infty} \sum_{m=n}^{(2k+1)^d} \mathbb{P}_\mu(|\mathcal{T}_k(0)| = m).$$

We will derive a bound for $\mathbb{P}_\mu(|\mathcal{T}_k(0)| = m)$. We write

$$\mathbb{P}_\mu(|\mathcal{T}_k(0)| = m) = \sum_{\substack{|\mathcal{C}|=m \\ 0 \in \mathcal{C}}} \mathbb{P}_\mu(\mathcal{T}_k(0) = \mathcal{C}),$$

where the sum runs over all finite connected subsets $\mathcal{C}$ of size $m$ containing the origin. Then, by Lemma 4.3, this implies a minimum number of at least $m - 1$ sand grains in $\mathcal{C}$ in $\eta_\infty$. But since no sand can have entered $\mathcal{C}$ during stabilization—in fact, grains must have left $\mathcal{C}$—it also implies that $\mathcal{C}$ contains at least $m$ grains at $t = 0$. Since $\mu$ is a product measure, this corresponds for $\rho < 1$ to a classical large deviation of $\frac{1}{m} \sum_{x \in \mathcal{C}} \eta(x)$, and we can bound the corresponding probability by a Chernov bound for sums of independent random variables, that is, there is a constant $\alpha = \alpha(\rho) > 0$ such that, for all $m$,

$$\mathbb{P}_\mu(\mathcal{T}_k(0) = \mathcal{C}) \leq e^{-\alpha m}. \tag{17}$$

(For this statement we do not yet need Lemma 4.1.) We now distinguish between (a) and (b):

(a) For $d = 1$, the number of clusters of size $m$ containing the origin is equal to $m$. Hence, for $d = 1$, we arrive at

$$\sum_{m=n}^{(2k+1)^d} \mathbb{P}_\mu(|\mathcal{T}_k(0)| = m) \leq \sum_{m=n}^{(2k+1)^d} \sum_{\substack{|\mathcal{C}|=m \\ 0 \in \mathcal{C}}} e^{-\alpha m} \leq \sum_{m=n}^{(2k+1)^d} m e^{-\alpha m} \leq e^{-c_1 n},$$

with $c_1$ positive for all $\rho < 1$. Since this outcome does not depend on $k$, when inserting this in (16), we obtain for $d = 1$

$$\mathbb{P}_\mu(|\mathcal{T}_\infty(0)| \geq n) \leq e^{-c_1 n},$$

proving (a).

(b) For $d > 1$, we have, according to Lemma 4.1, that $\lim_{\rho \downarrow 0} \alpha(\rho) = \infty$. Also, for $d > 1$, there is a constant $\alpha' = \alpha'(d)$ such that the number of clusters of size $m$ containing the origin is at most $e^{\alpha' m}$; see, for example, [5]. Hence, we calculate

$$\sum_{m=n}^{(2k+1)^d} \mathbb{P}_{\mu_\rho}(|\mathcal{T}_k(0)| = m) \leq \sum_{m=n}^{(2k+1)^d} \sum_{\substack{|\mathcal{C}|=m \\ 0 \in \mathcal{C}}} e^{-\alpha m} \leq \sum_{m=n}^{(2k+1)^d} e^{(\alpha' - \alpha) m} \leq e^{-c_d n},$$

with $c_d$ positive for $\rho$ small enough. The proof is now finished as in case (a). □

REMARK 4.4. In the proof of Theorem 4.2, we used that $\mu$ is translation invariant, and that we have, for $\rho$ small enough, a large deviation bound for sums like $\sum_{x \in \Lambda} \eta(x)$, with $\Lambda$ some connected volume in $\mathbb{Z}^d$. There are many more measures that satisfy these requirements, for instance, Gibbs measures or other sufficiently rapidly mixing measures.

The argument to prove the exponential tail of the distribution of $|\mathcal{W}_\infty(0)|$, which in turn implies the exponential tail of the distribution of $|\mathcal{V}_\infty(0)|$, is similar, although some extra arguments are needed.

THEOREM 4.5. *Let $\mu_\rho$, $0 < \rho < 1$ be a collection of product measures satisfying conditions* (i)–(iii) *of Lemma 4.1. Then for all $\rho$ sufficiently small, there exists a constant $\gamma_d = \gamma_d(d, \rho) > 0$ such that*

$$\mathbb{P}_{\mu_\rho}(|\mathcal{W}_\infty(0)| \geq n) \leq e^{-\gamma_d n}.$$

PROOF. As in the proof of Theorem 4.2, we stabilize $\eta$ in nested boxes $B_k$, and write [see (16)]

$$\mathbb{P}_\mu(|\mathcal{W}_\infty(0)| \geq n) = \lim_{k \to \infty} \mathbb{P}_\mu(|\mathcal{W}_k(0)| \geq n)$$
$$= \lim_{k \to \infty} \left( \sum_{m=n}^{\infty} \mathbb{P}_\mu(|\mathcal{W}_k(0)| = m) + \mathbb{P}_\mu(|\mathcal{W}_k(0)| = \infty) \right).$$

The cluster $\mathcal{W}_k(0)$ consists of the following types of sites: sites that have toppled, sites that did not topple but received at least one grain, and sites that did not topple nor received grains but which were nonempty in $\eta$. The first two types of sites we can only find in the box $B_{k+1}$, but the third type we can also find outside this box. Outside the box $B_{k+1}$, the configuration did not change yet, so restricted to $\mathbb{Z}^d \setminus B_{k+1}$, we just have independent site percolation of nonempty sites. We take $\rho$ so small that the density of nonempty sites is below the critical value for independent site percolation, so that for every $k$, $|\mathcal{W}_k(0)|$ is finite a.s. We write

$$\mathbb{P}_\mu(|\mathcal{W}_\infty(0)| \geq n) = \lim_{k \to \infty} \sum_{m=n}^{\infty} \mathbb{P}_\mu(|\mathcal{W}_k(0)| = m)$$
$$= \lim_{k \to \infty} \sum_{m=n}^{\infty} \sum_{\substack{|\mathcal{C}|=m \\ 0 \in \mathcal{C}}} \mathbb{P}_\mu(\mathcal{W}_k(0) = \mathcal{C}),$$

and again derive a bound for $\mathbb{P}_\mu(\mathcal{W}_k(0) = \mathcal{C})$ using that, on the event $\{\mathcal{W}_k(0) = \mathcal{C}\}$, there must have been a certain minimal number of sand grains in $\mathcal{C}$ before stabilization. Suppose $\mathcal{W}_k(0) = \mathcal{C}$. If $\mathcal{C}$ contains a cluster of size $m_t \geq 1$ of toppled sites, with $m_b \geq 2d$ boundary sites, then the number of grains in this region of sites—after as well as before toppling—is at least $m_t - 1 + m_b \geq 2d$, so that the density in this region is at least $\frac{2d}{2d+1}$. $\mathcal{C}$ might contain several of these regions, as well as nonempty sites that did not topple nor receive any grains. Thus, we cannot conclude more than that the density in $\mathcal{C}$ before toppling was at least $\frac{2d}{2d+1}$, which for $\rho < \frac{2d}{2d+1}$ corresponds to a large deviation of $\frac{1}{m} \sum_{x \in \mathcal{C}} \eta(x)$.

The rest of the proof proceeds the same as for Theorem 4.2. Note that the fact that we now sum $m$ from $n$ to $\infty$ instead of to $(2k+1)^d$ makes no difference for the outcome. $\square$

REMARK 4.6. For $d = 1$, the critical density of nonempty sites is 1, but for all $\rho$ we have that $\mathbb{P}_\mu(\eta(0) = 1) < 1$. Therefore, Theorem 4.5 is valid for $\rho < \frac{2}{3}$. However, in $d = 1$ it is not hard to see that, for all $\rho < 1$, $|\mathcal{W}_\infty(0)|$ is finite a.s. Indeed, it is not hard to see that there a.s. is a positive density of pairs of neighboring sites which never topple (compare the last part of Theorem 2.8). This implies that $|\mathcal{W}_\infty(0)|$ is finite.

**5. Some open problems.** It is clear that our results are a first step in the study of sandpile percolation and that many challenging open problems remain. We mention some of them:

1. *Infinite sandpile percolation clusters.* Starting from a product measure $\mu_\rho$ with density $\rho$ on $\mathcal{X}$, we know that, for $\rho > 2d - 1$, the measure is not stabilizable, and all sites will topple infinitely many times. For small $\rho$ we have sandpile percolation clusters that look like subcritical clusters (of ordinary percolation). Are there values of $\rho < 2d - 1$ such that there is an infinite cluster of toppled sites, but $\mu_\rho$ is still stabilizable? A guess would be that this happens as soon as $\rho > d$, the density of minimally recurrent configurations. If such a percolation transition occurs, is there a unique infinite cluster?
2. *Critical value.* Let $\mathcal{M}$ be a one-parameter family of Poisson product measures with parameter $\rho$. Is $\rho_c(\mathcal{M}) < 2d - 1$? Since it is conjectured that $\rho_c$ is the expected height in the critical sandpile model, in the limit of large volumes, this inequality should hold. However, we expect that, especially in high $d$, it will be easier to show a product measure with expected height "close to $2d - 1$," and not concentrating on stable configurations, is not stabilizable.
3. *Behavior at the critical value.* In $d = 1$ we proved that the Poisson measure $\mu_\rho$ is not stabilizable at the critical value $\rho = 1$. Is it true that, for all $d \geq 1$, $\mu_\rho$ is not stabilizable at the critical value $\rho = \rho_c$, where $\rho_c$ is the critical value of the family of Poisson measures?

A. Fey
R. Meester
VU—University Amsterdam
De Boelelaan 1081a
1081 HV Amsterdam
The Netherlands
E-mail: fey@eurandom.tue.nl
    rmeester@few.vu.nl

F. Redig
Mathematisch Instituut Universiteit Leiden
Niels Bohrweg 1
2333 CA Leiden
The Netherlands
E-mail: redig@math.leidenuniv.nl